\documentclass[12pt,reqno]{amsart}
\usepackage{amssymb}

 \rm

\renewcommand{\(}{\left(}
\renewcommand{\)}{\right)}
\renewcommand{\[}{\left[}
\renewcommand{\]}{\right]}
\newcommand{\<}{\langle}
\renewcommand{\>}{\rangle}

\newcommand{\norm}[1]{\left\lVert#1\right\rVert}

\renewcommand{\phi}{\varphi}

\renewcommand{\Re}{{\mathrm{Re}\,}}

\newcommand{\CC}{\mathbb{C}}
\newcommand{\RR}{\mathbb{R}}

\newcommand{\NN}{\mathbb{N}}

\newcommand{\Res}{\operatorname{Res}}

\theoremstyle{plain}
\newtheorem{thm}{Theorem}[section]
\newtheorem{lem}[thm]{Lemma}

\theoremstyle{definition}
\newtheorem{defn}[thm]{Definition}

\theoremstyle{remark}
\newtheorem{rem}[thm]{Remark}

\title{The Brylinski beta function of a coaxial layer}
\author{Pooja Rani and M.~K.~Vemuri}
\address{Department of Mathematical Sciences, IIT (BHU), Varanasi 221 005.}
\keywords{Analytic continuation; distribution; energies of submanifolds;
  invariant theory; second fundamental form.\\
  Pooja Rani is supported by DST-INSPIRE Fellowship/2019/IF190300.}

\begin{document}

\begin{abstract}
In [Pooja Rani and M.~K.~Vemuri, The Brylinski beta function of a
double layer, Differential Geom. Appl. \textbf{92}(2024)], an
analogue of Brylinski's knot beta function was defined for a compactly
supported (Schwartz) distribution $T$ on Euclidean space.  Here we
consider the Brylinski beta function of the distribution defined by a
coaxial layer on a submanifold of Euclidean space.  We prove that it
has an analytic continuation to the whole complex plane as a
meromorphic function with only simple poles, and in the case of a
coaxial layer on a space curves, we compute some of the residues in
terms of the curvature and torsion.

\end{abstract}

\maketitle

\section{Introduction}\label{S:intro}

In \cite{bry}, Brylinski introduced the beta function of a geometric knot in
$\RR^3$.  He was partly motivated by the desire to give a definition of
M\"obius energy (see \cite{fre,OHara1991241}) independent of an arbitrary
``renormalization''.  However, he also gave some beautiful formulae for the
first few residues of his beta function.  They turn out to be integrals of
polynomials in the curvature, torsion and their derivatives.

The basic idea behind Brylinski's beta function is much more general.  For
instance, in \cite{bbs}, Fuller and Vemuri defined and studied the beta
function of a compact submanifold of Euclidean space
(see also \cite{o2018regularized, err-o2018regularized}).
In \cite{bbdl}, Rani and Vemuri defined the beta function of a compactly
supported (Schwartz) distribution on Euclidean space, and studied
the beta function of a uniform double layer distribution on a
hypersurface.

Here we consider the beta function of a uniform coaxial layer on a
compact submanifold of Euclidean space.  In Section \ref{S:BBF}, we
define the concepts of coaxial derivative and coaxial layer,
and give an explicit expression for the beta function of a coaxial
layer.
In Section \ref{S:AC}, we express the beta function of a coaxial
layer in terms of Fuller and Vemuri's beta function and an additional
term involving the mean curvature.  We use this to give the analytic
continuation of the Brylinski beta function of a coaxial layer to the
whole complex plane as a
meromorphic function with only simple poles.
In Section \ref{S:curves}, we consider, in detail, the case of curves
in three dimensional space, and compute the residues at $3$, $1$ and $-1$.

\section{The Brylinski beta function}\label{S:BBF}

Let $T$ be a compactly supported distribution on $\RR^d$.  Since $T$
has finite order, there exists $k \in \NN$
such that $T$ extends to a continuous linear functional on $C^k(\RR^d)$.
Observe that if $\Re s > 2k$ and $N_s(u)=\norm{u}^s$ then
$N_s \in C^{2k}(\RR^d)$, and so $N_s*T \in C^k(\RR^d)$.
The following definition was given in \cite{bbdl}.

\begin{defn}
The Brylinski beta function of $T$ is the function $B_T(s)$, defined for
$\Re s > 2k$ by
\begin{equation*}
B_T(s) = T(N_s * T).
\end{equation*}
\end{defn}

Suppose $M \subseteq \RR^d$ is a compact smooth $n$-dimensional
submanifold, and $dA$ is the surface measure on $M$.

If $T:C_0^\infty(\RR^d)\to\CC$ is the distribution defined by
\begin{equation*}
T\phi=\int_M \phi \, dA,
\end{equation*}
then $B_T=B_M$ is the beta function studied by Fuller and Vemuri in \cite{bbs}.

Given $u\in M$, let $\nu_u=(T_uM)^\perp$ be the normal space to $M$ at $u$.
Let $\Sigma_r$ denote the sphere of radius $r$ centered at the origin in
$\nu_u$.

\begin{lem}
If $\phi$ is a $C^2$ function defined in a neighborhood of $u$, then
\begin{equation*}
\lim_{r\to 0} 
\frac{\frac{1}{\omega_{d-n-1} r^{d-n-1}}
\int_{\Sigma_r} \phi(u+\sigma) \, d\sigma - \phi(u)}{r^2}
\end{equation*}
exists.
\end{lem}

\begin{proof}
Let $X_1, \dots, X_{d-n}$ be orthonormal vectors which span $\nu_u$, and
extend them to constant vector fields on $\RR^d$.  We claim that the limit
is in fact
\begin{equation*}
\frac{(X_1^2\phi + \cdots +X_{d-n}^2\phi)(u)}{2(d-n)}.
\end{equation*}

By Taylor's theorem, we may write, for $\sigma \in \Sigma_r$,
\begin{equation*}
\phi(u+\sigma)=
\phi(u)+\sum_{j=1}^{d-n} (X_j\phi)(u)\, \<\sigma, X_j\>
+\frac12 \sum_{j,k=1}^{d-n} (X_jX_k\phi)(u) \, \<\sigma, X_j\> \, \<\sigma, X_k\>
+o(r^2).
\end{equation*}

Therefore
\begin{equation*}
\begin{gathered}
\frac{1}{\omega_{d-n-1} r^{d-n-1}}
\int_{\Sigma_r} \phi(u+\sigma) \, d\sigma - \phi(u)\\
=
\frac{1}{\omega_{d-n-1} r^{d-n-1}}
\(
\begin{gathered}
\sum_{j=1}^{d-n} (X_j\phi)(u)\, \int_{\Sigma_r} \<\sigma, X_j\> \, d\sigma\\
+\frac12 \sum_{j,k=1}^{d-n} (X_jX_k\phi)(u) \,
         \int_{\Sigma_r}\<\sigma, X_j\> \, \<\sigma, X_k\> \, d\sigma\\
\end{gathered}
\)
+ o(r^2)\\
=
\frac{(X_1^2\phi + \cdots +X_{d-n}^2\phi)(u)}{2(d-n)} r^2 + o(r^2),
\end{gathered}
\end{equation*}
because
\begin{equation*}
\begin{aligned}
\int_{\Sigma_r} \<\sigma, X_j\> \, d\sigma
=&\; 0, \quad\text{and}\\
\frac{1}{\omega_{d-n-1}r^{d-n+1}}
\int_{\Sigma_r}\<\sigma, X_j\> \, \<\sigma, X_k\> \, d\sigma
=&\; \frac{\delta_{jk}}{d-n}.
\end{aligned}
\end{equation*}
\end{proof}

We define the {\em coaxial derivative} of $\phi$ at $u$ by
\begin{equation*}
\Delta_{\circ_u}\phi = 2(d-n)\lim_{r\to 0} 
\frac{\frac{1}{\omega_{d-n-1} r^{d-n-1}}
\int_{\Sigma_r} \phi(u+\sigma) \, d\sigma - \phi(u)}{r^2}.
\end{equation*}

Thus, $\Delta_{\circ_u}=X_1^2+\cdots+X_{d-n}^2$
where $\{X_1, \dots,
X_{d-n}\}$ is any orthonormal basis of $\nu_u$, and so
$\Delta_{\circ_u}$ is a distribution of order two supported at $u$.
Define $T_\circ:C_0^\infty \to \CC$ by
\begin{equation*}
T_\circ\phi=\int_M \Delta_{\circ_u}\phi \, dA(u).
\end{equation*}
Then $T_\circ$ is a distribution of order two supported on $M$, which we call
the {\em coaxial layer} on $M$.

If $\Re s > 4$, then
\begin{equation*}
(N_s * T_\circ)(u) = T_\circ(\tau_{-u}N_s) = \int_M \Delta_{\circ_v}N_s(v-u) \, dA(v),
\end{equation*}
and so
\begin{equation*}
B_{T_\circ}(s) = T_\circ(N_s * T_\circ)
= \int_M \int_M \Delta_{\circ_u}\Delta_{\circ_v}N_s(v-u)
\, dA(v) \, dA(u).
\end{equation*}
We will refer to $B_{T_\circ}$ as the
{\em Brylinski beta function of a coaxial layer on $M$}, and will use the
notation
\begin{equation*}
B_{M_\circ}=B_{T_\circ}.
\end{equation*}

\begin{rem}
In calculations, it is convenient to regard $N_s$ as a function of two
variables:
\begin{equation*}
N_s(u,v)=N_s(u-v).
\end{equation*}
\end{rem}

\section{The analytic continuation of $B_{M_\circ}$}\label{S:AC}

\begin{thm}\label{T:AC}
The function $B_{M_\circ}$ admits an analytic continuation to the whole complex
plane as a meromorphic function with only simple poles at
$-n-j$, $j=-4,-2,0,2,4,\dots$.
\end{thm}

\begin{proof}
Observe first that if $\phi\in C_0^\infty(\RR^d)$ then
\begin{equation*}
\Delta_{\circ_u} \phi = (\Delta\phi - \Delta_M\phi + H\phi)(u),
\end{equation*}
(see\cite{coaxial}) where $\Delta$ is the Laplacian on $\RR^d$, $\Delta_M$ is the
Laplace-Beltrami operator of the Riemannian manifold $M$
(see \cite[p99]{chavel}) and $H$ is the mean curvature vector of $M$.
Therefore
\begin{equation*}
\begin{aligned}
T_\circ\phi
=&\; \int_M \Delta\phi \, dA(u) + \int_M H\phi \, dA(u),
\end{aligned}
\end{equation*}
because
\begin{equation*}
\int_M \Delta_M \phi \, dA(u) = 0,
\end{equation*}
by the divergence theorem (see \cite[p149]{chavel}).
Therefore, if $\Re s > 4$, then
\begin{equation*}
(N_s * T_\circ)(u) = T_\circ(\tau_{-u}N_s) = \int_M \Delta_v N_s(v-u) \, dA(v)
+ \int_M H_v N_s(v-u) \, dA(v),
\end{equation*}
and so
\begin{equation}\label{{E:bbcl}}
\begin{aligned}
B_{T_\circ}(s)
=&\; T_\circ(N_s * T_\circ)\\
=&\; \int_M \int_M \Delta_u\Delta_v N_s(v-u) \, dA(v) \, dA(u)+
\int_M \int_M H_uH_v N_s(v-u) \, dA(v) \, dA(u)\\
=&\; B_1(s) + B_2(s),
\end{aligned}
\end{equation}
where
\begin{equation*}
\begin{aligned}
B_1(s) =&\; \int_M \int_M \Delta_u\Delta_v N_s(v-u) \, dA(v) \, dA(u),
\qquad\text{and}\\
B_2(s) =&\; \int_M \int_M H_uH_v N_s(v-u) \, dA(v) \, dA(u),
\end{aligned}
\end{equation*}
because
\begin{equation*}
\Delta_vH_u N_s(v-u) + H_v\Delta_u N_s(v-u) = 0.
\end{equation*}

For $v\in M$, let
\begin{equation*}
\begin{aligned}
B_M^v =&\; \int_M N_s(v-u) \, dA(u),\\
B_1^v =&\; \int_M \Delta_u\Delta_v N_s(v-u) \, dA(u), \qquad\text{and}\\
B_2^v =&\; \int_M H_u H_v N_s(v-u) \, dA(u).
\end{aligned}
\end{equation*}
By a translation of $M$, we may assume that $v=0$, and so
\begin{equation*}
\begin{aligned}
B_1^v
=&\; s(s-2+d)(s-2)(s-4+d) \int_M \(\sum_{i=1}^d u_i^2\)^{\frac{s-4}{2}} \, dA(u)\\
=&\; s(s-2+d)(s-2)(s-4+d) B_M^v(s-4).
\end{aligned}
\end{equation*}
Integration over $M$ with respect to $v$ gives
\begin{equation*}
B_1(s) = s(s-2)(s-2+d)(s-4+d)B_M(s-4), \qquad \Re s > 4.
\end{equation*}
It follows from \cite[Theorem 3.3]{bbs}, that $B_1$ has an analytic
continuation to the whole complex plane as a meromorphic function with
possible simple poles at $4-n-k$, $k=0, 2, 4, \dots$.

We claim that
$B_2$ may be analytically continued to a meromorphic function on $\CC$ with
simple poles at $-n-j$, $j=-2, 0, 2, 4, \dots$.  It suffices to give the
argument for $B_2^v$ (integration over $M$ gives the result for $B_2$).

Let $\psi\in C_c^\infty(\RR^d)$ be identically $1$ on a neighborhood of $v$, and
define
\begin{equation*}
B_2^\psi(s) = \int_M H_u H_v N_s(u,v) \psi(u) \, dA(u).
\end{equation*}
Then $B_2^\psi$ has the same principal part as $B_2^v$, so it suffices to
prove our claim with $B_2^v$ replaced by $B_2^\psi$.

By rotating $M$, we may assume that the tangent space to $M$ at $v=0$ is
$\RR^n \subseteq \RR^d$ (clearly this process does not affect $B_2$).
Then, in a neighborhood of $0$, $M$ is the graph of a function
$f:\RR^n \to \RR^{d-n}$ which vanishes to second order at $0$.  By
making the neighborhood smaller, we may assume that
$\norm{f(w)} < \norm{w}$.  Choose $\psi$ to have support in that neighborhood.
Let $\phi(w)=\psi(w, f(w)) A(w)$, where $A(w)$ is the area-density of $M$ (it
may be expressed in terms of the partial derivatives of $f$).

Write
\begin{equation*}
H_u=a_1\frac{\partial}{\partial u_1} + \cdots + a_d\frac{\partial}{\partial u_d},
\end{equation*}
and
\begin{equation*}
H_v=a_1\frac{\partial}{\partial v_1} + \cdots + a_d\frac{\partial}{\partial v_d},
\end{equation*}
where $a_1, \dots, a_d$ are smooth functions defined on $M$.

Then
\begin{equation*}
\begin{aligned}
H_u H_v N_s(u, 0)
=&\; -s\(\sum_{i=1}^d a_i(0)a_i(u)\)\(\sum_{i=1}^d u_i^2\)^{\frac{s-2}{2}}\\
 &\; -s(s-2)\(\sum_{i=1}^d a_i(0) u_i\)\(\sum_{i=1}^d a_i(u)u_i\)
            \(\sum_{i=1}^d u_i^2\)^{\frac{s-4}{2}}\\
=&\; -s\(H_0\cdot H_{(w,f(w))}\)\(\norm{w}^2 + \norm{f(w)}^2\)^{\frac{s-2}{2}}\\
 &\; -s(s-2)\(H_0\cdot(w,f(w))\)\(H_{(w,f(w))}\cdot(w,f(w))\)
\(\norm{w}^2 + \norm{f(w)^2}\)^{\frac{s-4}{2}},
\end{aligned}
\end{equation*}
where $X\cdot Y$ denotes the dot product of $X$ and $Y$ in $\RR^d$.\\
Therefore
\begin{equation*}
\begin{aligned}
B_2^\psi(s)  
=&\;  -s(s-2)\int_{\RR^n}\norm{w}^{s}
             \(1+{\frac{f(w)^2}{\norm{w}^2}}\)^{\frac{s-4}{2}}
             \frac{\(H_0\cdot(w,f(w))\)}{\norm{w}^2}
             \frac{\(H_{(w,f(w))}\cdot(w,f(w))\)}{\norm{w}^2}
             \,\phi(w)\, dw\\
 &\;   -s\int_{\RR^n}\norm{w}^{s-2}
         \(1+{\frac{f(w)^2}{\norm{w}^2}}\)^{\frac{s-2}{2}} 
         \(H_0\cdot H_{(w,f(w))}\)\,\phi(w)\, dw\\
=&\; -s(s-2)\int_0^\infty r^{s+n-1} S_1(r,s)\,dr 
     - s\int_0^\infty r^{s+n-3} S_2(r,s)\,dr,
\end{aligned}
\end{equation*}
where
\begin{equation*}
\begin{gathered}
S_1(r,s)=
 \int_{S^{n-1}(1)}\(1+{\frac{f(rw)^2}{r^2}}\)^{\frac{s-4}{2}}
 \frac{\(H_0\cdot(w,f(w))\)}{\norm{w}^2}
  \frac{\(H_{(w,f(w))}\cdot(w,f(w))\)}{\norm{w}^2}
\,\phi(w)\, d\sigma(w),\\
S_2(r,s)=
 \int_{S^{n-1}(1)}\(1+{\frac{f(rw)^2}{r^2}}\)^{\frac{s-2}{2}}  
 \(H_0\cdot H_{(w,f(w))}\) \,\phi(w)\, d\sigma(w),
\end{gathered}
\end{equation*}
and $\sigma$ is the surface measure on $S^{n-1}(1)$.
Note that
$\(H_0\cdot(w,f(w))\)$ and $\(H_{(w,f(w))}\cdot(w,f(w))\)$ vanish to
order two at $0$ because $H$ is orthogonal to $M$.

Using the same argument as in the proof of \cite[Theorem 3.1]{bbdl},
it follows that $B_2^v$ may be analytically continued to a meromorphic
function on $\CC$ with simple poles at $-n-j$, $j=-2, 0, 2, 4, \dots$,
and in particular, the residues of $B_2$ may be computed by the formulae
\begin{equation*}
\begin{gathered}
 \Res_{s=2-n} B_2^v = -sS_2(0,s),
 \qquad\text{and}\\
\Res_{s=-n-j} B_2^v
= -s(s-2)\frac{1}{j!} \frac{\partial^j S_1}{\partial r^j}(0, -n-j)
  -s \frac{1}{j+2!} \frac{\partial^{j+2} S_2}{\partial r^{j+2}}(0, -n-j).
\end{gathered}
\end{equation*}

The result now follows from Equation (\ref{{E:bbcl}}).
\end{proof}

\section{Curves in $\RR^3$}\label{S:curves}

In this section, we assume $M \subseteq \RR^3$ is a smooth curve.  We will
compute $\mathrm{Res}_{s=k} B_{M_\circ}$ for $k=3$, $1$, and $-1$ in terms
of the curvature, torsion and their derivatives with respect to
a unit speed parametrization. 

\begin{thm}
\begin{equation*}
\begin{aligned}
\Res_{s=3} B_{M_\circ} &= 48 \, \mathrm{Length}(M)\\    
\Res_{s=1} B_{M_\circ} &= -2\int_M  \kappa_0^2  \, dA\\
\Res_{s=-1} B_{M_\circ} &= \int_M \(
\frac{3}{4} \kappa_0^4 - \kappa_0^2 \tau_0^2 
+ \kappa_0 \kappa_2 
\) \, dA
\end{aligned}
\end{equation*}
where $\kappa_0$ and $\tau_0$ are the curvature and torsion of $M$ respectively,
and $\kappa_n$ and $\tau_n$ are the $n$-th
derivative of $\kappa_0$ and $\tau_0$ with respect to a
unit speed parametrization.
\end{thm}

\begin{proof}
As in the proof of Theorem \ref{T:AC}, we compute the residues of
$B_{M_\circ}^v$ for each $v \in M$; integration over $M$ gives the stated
result. Also, we may assume that $v=0$, and locally $M$ is the graph of a
function $f=(f_1,f_2)$ which vanishes to
second order at $0$.  By the proof of Theorem \ref{T:AC},
\begin{equation*}
B_1^v(s) = (s+1)s(s-1)(s-2)B_M^v(s-4),
\end{equation*}
and by the proof of Theorem 3.3 of \cite{bbs},
\begin{equation*}
\Res_{s=-1-j} B_M^v
=\frac{1}{j!} \frac{\partial^j S}{\partial r^j}(0, -1-j),
\end{equation*}
where
\begin{equation*}
S(r,s)=\sum_{w=\pm 1}\(1+{\frac{f(rw)^2}{r^2}}\)^{\frac{s}{2}}  \,\phi(rw),\\   
\end{equation*}
for small $r$.

If we write
\begin{equation}\label{E:monge}
\begin{gathered}
f_1(u_1) =
a_2 u_1^2 + a_3 u_1^3 + a_4 u_1^4 +
a_5 u_1^5 + a_6 u_1^6 + a_7 u_1^7 + a_8 u_1^8 
+ O(u_1^9),\\
f_2(u_1) =
b_2 u_1^2 + b_3 u_1^3 + b_4 u_1^4 +
b_5 u_1^5 + b_6 u_1^6 + b_7 u_1^7 + b_8 u_1^8 
+ O(u_1^9),\\
\(1+\frac{f(rw)^2}{r^2}\)^{\frac{s}{2}}
=
\(1 + \frac{s}{2} \frac{f^2}{r^2} + \frac{s(s-2)}{8} \frac{f^4}{r^4}+
\frac{s(s-2)(s-4)}{48} \frac{f^6}{r^6}\)
+ O(r^8),
\end{gathered}
\end{equation}
substitute into the definition of $S(r,s)$, and
perform the indicated differentiation (using the computer algebra package
Maxima \cite{maxima}), we find
\begin{equation*}
\begin{aligned}
\Res_{s=-1} B_M^v
=&\; 2,  \\
\Res_{s=-3} B_M^v
=&\; b_{2}^2+a_{2}^2,  \quad\text{and}\\
\Res_{s=-5} B_M^v
=&\; {{24\,b_{2}\,b_{4}+16\,b_{3}^2-21\,b_{2}^4-42\,a_{2}^2\,b_{2}^2+24
     \,a_{2}\,a_{4}+16\,a_{3}^2-21\,a_{2}^4}\over{4}}.\\
\end{aligned}
\end{equation*}

Therefore,
\begin{equation*}
\begin{aligned}
\Res_{s=3} B_1^v
=&\; 48,\\
\Res_{s=1} B_1^v
=&\; 0, \quad\text{and} \\
\Res_{s=-1} B_1^v
=&\; 0.\\
\end{aligned}
\end{equation*}
Now,
\begin{equation*}
\begin{gathered}     
\Res_{s=1} B_2^v = -sS_2(0,s),
\qquad\text{and}\\
\Res_{s=-1-j} B_2^v
=-s(s-2)\frac{1}{j!} \frac{\partial^j S_1}{\partial r^j}(0, -1-j)
 -s\frac{1}{j+2!} \frac{\partial^{j+2} S_2}{\partial r^{j+2}}(0, -1-j),
\end{gathered}
\end{equation*}
where
\begin{equation*}
\begin{gathered}
S_1(r,s)
=\sum_{w=\pm 1}\(1+{\frac{f(rw)^2}{r^2}}\)^{\frac{s-4}{2}}
\frac{\(H_0\cdot(rw,f(rw))\)}{\norm{w}^2}
\frac{\(H_{(rw,f(rw))}\cdot(rw,f(rw))\)}{\norm{w}^2}
\,\phi(rw),
\qquad\text{and}\\
S_2(r,s)=\sum_{w=\pm 1}\(1+{\frac{f(rw)^2}{r^2}}\)^{\frac{s-2}{2}}
\(H_0\cdot H_{(rw,f(rw))}\)\,\phi(rw),
\end{gathered}
\end{equation*}
for small $r$,
with
\begin{equation*}
H_{(u_1,f(u_1))} =
\frac{\(  
-f_1'f_1''-f_2'f_2'',\,
f_2'^2f_1''-f_1'f_2'f_2''+
f_1'',\,
f_1'^2f_2''-f_2'f_1'f_1''+
f_2''
\)
}{(1+f_1'^2+f_2'^2)^2}. 
\end{equation*}

If we write $f_1$ and $f_2$ as indicated in equation (\ref{E:monge}),
\begin{equation*}
\begin{gathered}
\(1+{\frac{f(rw)^2}{r^2}}\)^{\frac{s-4}{2}}
\frac{\(H_0\cdot(rw,f(rw))\)}{r^2}
\frac{\(H_{(rw,f(rw))}\cdot(rw,f(rw))\)}{r^2}\\
=\(1 + \frac{(s-4)}{2} \frac{f^2}{r^2} + \frac{(s-4)(s-6)}{8} \frac{f^4}{r^4}+
\frac{(s-4)(s-6)(s-8)}{48} \frac{f^6}{r^6}\)\frac{\(H_0\cdot(rw,f(rw))\)}{r^2}\\
\frac{\(H_{(rw,f(rw))}\cdot(rw,f(rw))\)}{r^2}+ O(r^{12}),
\end{gathered}
\end{equation*}
and
\begin{equation*}
\begin{gathered}
\(1+{\frac{f(rw)^2}{r^2}}\)^{\frac{s-2}{2}}\(H_0\cdot H_{(rw,f(rw))}\)\\
=\(1 +  \frac{(s-2)}{2} \frac{f^2}{r^2} + \frac{(s-2)(s-4)}{8} \frac{f^4}{r^4} +
\frac{(s-2)(s-4)(s-6)}{48} \frac{f^6}{r^6}\)\(H_0\cdot H_{(rw,f(rw))}\)+ O(r^{10}),
\end{gathered}
\end{equation*}
substitute into the definition of $S_1(r,s)$ and $S_2(r,s)$ and
perform the indicated differentiation (using the computer algebra package
Maxima \cite{maxima}), we find
\begin{equation*}
\begin{aligned}
\Res_{s=1} B_2^v
=&\; -8\,b_{2}^2-8\,a_{2}^2,  \quad\text{and}\\
\Res_{s=-1} B_2^v
=&\; 48\,b_{2}\,b_{4}-36\,b_{2}^4-72\,a_{2}^2\,b_{2}^2+48\,a_{2}\,a_{4}-
 36\,a_{2}^4.\\
\end{aligned}
\end{equation*}
Therefore
\begin{equation*}
\begin{aligned}
\Res_{s=3} B_{M_\circ}^v
=&\; 48,\\
\Res_{s=1} B_{M_\circ}^v
=&\; -8\,b_{2}^2-8\,a_{2}^2, \quad\text{and}\\
\Res_{s=-1} B_{M_\circ}^v
=&\; 48\,b_{2}\,b_{4}-36\,b_{2}^4-72\,a_{2}^2\,b_{2}^2+48\,a_{2}\,a_{4}-
 36\,a_{2}^4.\\
\end{aligned}
\end{equation*}

At $u_1=0$, we compute (using Maxima)
\begin{equation*}
\begin{aligned}
\kappa_0
=&\; \sqrt{4\,b_{2}^2+4\,a_{2}^2}, \\
\kappa_1
=&\; {{24\,b_{2}\,b_{3}+24\,a_{2}\,a_{3}}\over{2\,\sqrt{4\,b_{2}^2+4\,
      a_{2}^2}}}, \\
\kappa_2
=&\; \frac{\sqrt{4\,b_{2}^2+4\,a_{2}^2}}
          {(b_{2}^2+a_{2}^2)^2}
\[
\begin{gathered}
( (12\,b_{2}^3+12\,a_{2}^2\,b_{2} )\,b_{4}+9\,a_{2}^2\,b_{3}^2
-18\,a_{2}\,a_{3}\,b_{2}\,b_{3}-12\,b_{2}^6-36\,a_{2}^2\,b_{2}^4\\
+ (12\,a_{2}\,a_{4}+9\,a_{3}^2-36\,a_{2}^4 )\,b_{2}^2
+12\,a_{2}^3\,a_{4}-12\,a_{2}^6)
\end{gathered}
\],\\
\kappa_3
=&\; \frac{\sqrt{4\,b_{2}^2+4\,a_{2}^2}}
          {(b_{2}^2+a_{2}^2)^3}
\[
\begin{gathered}
( (60\,b_{2}^5+120\,a_{2}^2\,b_{2}^3+60\,a_{2}^4\,b_{2} )\,b_{5}
+ ( (108\,a_{2}^2\,b_{2}^2+108\,a_{2}^4 )\,b_{3}-108\,a_{2}\,a_{3}\,b_{2}^3\\
-108\,a_{2}^3\,a_{3}\,b_{2} )\,b_{4}-81\,a_{2}^2\,b_{2}\,b_{3}^3+
(162\,a_{2}\,a_{3}\,b_{2}^2-81\,a_{2}^3\,a_{3} )\,b_{3}^2+
(-228\,b_{2}^7-684\,a_{2}^2\,b_{2}^5\\
+ (-108\,a_{2}\,a_{4}-81\,a_{3}^2-684\,a_{2}^4 )\,b_{2}^3
+ (-108\,a_{2}^3\,a_{4}+162\,a_{2}^2\,a_{3}^2-228\,a_{2}^6 )\,b_{2} )\,b_{3}\\
-228\,a_{2}\,a_{3}\,b_{2}^6+ (60\,a_{2}\,a_{5}+108\,a_{3}\,a_{4}-684\,
a_{2}^3\,a_{3} )\,b_{2}^4+ (120\,a_{2}^3\,a_{5}+108\,a_{2}^2\,a_{3}\,a_{4}\\
-81\,a_{2}\,a_{3}^3-684\,a_{2}^5\,a_{3} )\,b_{2}^2
+60\,a_{2}^5\,a_{5}-228\,a_{2}^7\,a_{3} )  
\end{gathered}
\],\\
\tau_0
=&\; {{3\,a_{2}\,b_{3}-3\,a_{3}\,b_{2}}\over{b_{2}^2+a_{2}^2}}, \\
\tau_1
=&\; \frac{1}{ (b_{2}^2+a_{2}^2)^2}
\[
\begin{gathered}
(12\,a_{2}\,b_{2}^2+12\,a_{2}^3 )\,b_{4}-18\,a_{2}\,
b_{2}\,b_{3}^2+ (18\,a_{3}\,b_{2}^2-18\,a_{2}^2\,a_{3} )\,b_{3}\\
-12\,a_{4}\,b_{2}^3+ (18\,a_{2}\,a_{3}^2-12\,a_{2}^2\,a_{4})\,b_{2}
\end{gathered}
\], \\
\tau_2
=&\; \frac{1}{ (b_{2}^2+a_{2}^2)^3}
\[
\begin{gathered}
(60\,a_{2}\,b_{2}^4+120\,a_{2}^3\,b_{2}^2+60\,a_{2}^5 )
\,b_{5}+ ( (-216\,a_{2}\,b_{2}^3-216\,a_{2}^3\,b_{2} )
\,b_{3}+108\,a_{3}\,b_{2}^4\\
-108\,a_{2}^4\,a_{3} )\,b_{4}+ (162\,a_{2}\,b_{2}^2-54\,a_{2}^3 )\,b_{3}^3
+ (486\,a_{2}^2\,a_{3}\,b_{2}-162\,a_{3}\,b_{2}^3 )\,b_{3}^2+ (96\,a_{2}\,b_{2}^6\\
+ (108\,a_{4}+288\,a_{2}^3 )\,b_{2}^4+ (288\,a_{2}^5-486\,a_{2}\,a_{3}^2 )\,b_{2}^2
-108\,a_{2}^4\,a_{4}+162\,a_{2}^3\,a_{3}^2\\
+96\,a_{2}^7 )\,b_{3}-96\,a_{3}\,b_{2}^7+(-60\,a_{5}-288\,a_{2}^2\,a_{3} )\,b_{2}^5
+ (-120\,a_{2}^2\,a_{5}+216\,a_{2}\,a_{3}\,a_{4}\\
+54\,a_{3}^3-288\,a_{2}^4\,a_{3} )\,b_{2}^3+ (-60\,a_{2}^4\,a_{5}
+216\,a_{2}^3\,a_{3}\,a_{4}-162\,a_{2}^2\,a_{3}^3-96\,a_{2}^6\,a_{3} )\,b_{2}
\end{gathered}
\],\\
\noalign{and}
\tau_3
=&\; \frac{1}{(b_{2}^2+a_{2}^2)^4}
\[
\begin{gathered}
(360\,a_{2}\,b_{2}^6+1080\,a_{2}^3\,b_{2}^4+1080\,a_{2}^5\,
b_{2}^2+360\,a_{2}^7 )\,b_{6}+ ( (-1440\,a_{2}\,b_{2}^5
-2880\,a_{2}^3\,b_{2}^3\\
-1440\,a_{2}^5\,b_{2} )\,b_{3}+720\,a_{3}\,b_{2}^6
+720\,a_{2}^2\,a_{3}\,b_{2}^4-720\,a_{2}^4\,a_{3}\,
b_{2}^2-720\,a_{2}^6\,a_{3} )\,b_{5}+ (-864\,a_{2}\,b_{2}^5\\
-1728\,a_{2}^3\,b_{2}^3-864\,a_{2}^5\,b_{2} )\,b_{4}^2+ ((3888\,a_{2}\,b_{2}^4
+2592\,a_{2}^3\,b_{2}^2-1296\,a_{2}^5 )\,b_{3}^2+ (-2592\,a_{3}\,b_{2}^5\\
+5184\,a_{2}^2\,a_{3}\,b_{2}^3+7776\,a_{2}^4\,a_{3}\,b_{2} )\,b_{3}
+1104\,a_{2}\,b_{2}^8+ (864\,a_{4}+4416\,a_{2}^3 )\,b_{2}^6
+ (864\,a_{2}^2\,a_{4}\\
-3888\,a_{2}\,a_{3}^2+6624\,a_{2}^5 )\,b_{2}^4+ (-864\,a_{2}^4\,a_{4}
-2592\,a_{2}^3\,a_{3}^2+4416\,a_{2}^7 )\,b_{2}^2-864\,a_{2}^6\,a_{4}\\
+1296\,a_{2}^5\,a_{3}^2+1104\,a_{2}^9)\,b_{4}+ (1944\,a_{2}^3\,b_{2}
-1944\,a_{2}\,b_{2}^3)\,b_{3}^4+ (1944\,a_{3}\,b_{2}^4-11664\,a_{2}^2\,a_{3}\,
b_{2}^2\\+1944\,a_{2}^4\,a_{3} )\,b_{3}^3+ (-792\,a_{2}\,
b_{2}^7+ (-1296\,a_{4}-2376\,a_{2}^3 )\,b_{2}^5+ (2592
\,a_{2}^2\,a_{4}+11664\,a_{2}\,a_{3}^2\\
-2376\,a_{2}^5 )\,b_{2}^3+ (3888\,a_{2}^4\,a_{4}-11664\,a_{2}^3\,a_{3}^2
-792\,a_{2}^7)\,b_{2} )\,b_{3}^2+ (792\,a_{3}\,b_{2}^8+ (720\,a_{5}\\
+1584\,a_{2}^2\,a_{3} )\,b_{2}^6+ (720\,a_{2}^2\,a_{5}
-7776\,a_{2}\,a_{3}\,a_{4}-1944\,a_{3}^3 )\,b_{2}^4+ (-720
\,a_{2}^4\,a_{5}-5184\,a_{2}^3\,a_{3}\,a_{4}\\
+11664\,a_{2}^2\,a_{3}^3-1584\,a_{2}^6\,a_{3} )\,b_{2}^2
-720\,a_{2}^6\,a_{5}+2592\,a_{2}^5\,a_{3}\,a_{4}-1944\,a_{2}^4\,a_{3}^3
-792\,a_{2}^8\,a_{3})\,b_{3}\\
-1104\,a_{4}\,b_{2}^9+ (-360\,a_{6}-4416\,a_{2}^2\,a_{4}
+792\,a_{2}\,a_{3}^2 )\,b_{2}^7+ (-1080\,a_{2}^2\,a_{6}+1440\,a_{2}\,a_{3}\,a_{5}\\
+864\,a_{2}\,a_{4}^2+ (1296\,a_{3}^2-6624\,a_{2}^4 )\,a_{4}
+2376\,a_{2}^3\,a_{3}^2 )\,b_{2}^5+ (-1080\,a_{2}^4\,a_{6}
+2880\,a_{2}^3\,a_{3}\,a_{5}\\
+1728\,a_{2}^3\,a_{4}^2+ (-2592\,a_{2}^2\,a_{3}^2-4416\,a_{2}^6)\,a_{4}
-1944\,a_{2}\,a_{3}^4+2376\,a_{2}^5\,a_{3}^2 )\,b_{2}^3+ (-360\,a_{2}^6\,a_{6}\\
+1440\,a_{2}^5\,a_{3}\,a_{5}+864\,a_{2}^5\,a_{4}^2+
(-3888\,a_{2}^4\,a_{3}^2-1104\,a_{2}^8 )\,a_{4}
+1944\,a_{2}^3\,a_{3}^4+792\,a_{2}^7\,a_{3}^2 )\,b_{2}
\end{gathered}
\].
\end{aligned}  
\end{equation*}

In order to re-express the residues in terms of the
curvature, torsion and their derivatives,
we use the concept of {\em weight} (see \cite[Theorem 5.1]{bbdl}).

Let us consider the weights of the various functionals at hand.

\begin{enumerate}

\item
The coefficients $a_i$ and $b_i$ have weight $i-1$.

\item
The functional $M\mapsto\Res_{s=-k} B_{M_\circ}^v$
has weight $k+3$.

\item
The invariants ${\kappa_n}(0)$ and ${\tau_n}(0)$
(in the $u_1$-coordinate system) have weight $n+1$.
\end{enumerate}

The following table lists the monomials in the invariants according to
weight.

\begin{center}
  \begin{tabular}{||c|c||}
    \hline
    \hline
    Weight & Monomial\\
    \hline
    $0$    & $1$\\
    \hline
    $1$    & $\kappa_0, \tau_0$\\
    \hline
    $2$    & $\kappa_0^2, \kappa_1, \tau_0^2, \tau_1, \kappa_0\tau_0$\\
    \hline
    $4$    & $\kappa_0^4, \kappa_1^2, \kappa_0\kappa_2, \kappa_0^2\kappa_1,
              \kappa_3, \tau_0^4, \tau_1^2, \tau_0\tau_2, \tau_0^2\tau_1$\\
           & $\tau_2\kappa_0, \tau_1\kappa_0^2, \tau_1\kappa_1, \kappa_2\tau_0,
              \kappa_1\tau_0^2, \kappa_0^2\tau_0^2, \tau_0\kappa_0^3,
              \kappa_0\tau_0^3,\tau_3$\\
       
    \hline
    \hline
  \end{tabular}
\end{center}

From this we conclude that $\Res_{s=3} B_{M_\circ}^v$ is independent
of $M$, $\Res_{s=1} B_{M_\circ}^v$ is a linear combination of
$\kappa_0^2$, $\kappa_1$, $\tau_0^2$, $\tau_1$, $\kappa_0\tau_0$,
and $\Res_{s=-1} B_{M_\circ}^v$ is a linear combination of
$\kappa_0^4$, $\kappa_1^2$, $\kappa_0\kappa_2$, $\kappa_0^2\kappa_1$,
$\kappa_3$, $\tau_0^4$, $\tau_1^2$, $\tau_0\tau_2$,
$\tau_0^2\tau_1$, $\tau_2\kappa_0$, $\tau_1\kappa_0^2$, $\tau_1\kappa_1$,
$\kappa_2\tau_0$, $\kappa_1\tau_0^2$, $\kappa_0^2\tau_0^2$, $\tau_0\kappa_0^3$,
$\kappa_0\tau_0^3$,$\tau_3$.

Using Macaulay2 \cite{M2}, we can find the coefficients:
\begin{equation*}
\begin{aligned}
\Res_{s=-1} B_M &= 2 \, \mathrm{Length}(M),\\    
\Res_{s=-3} B_M &= \int_M \frac{1}{4}\kappa_0^2  \, dA, \quad\text{and}\\
\Res_{s=-5} B_M &= \int_M \(
\frac{3}{64} \kappa_0^4 - \frac{1}{72}\kappa_0^2 \tau_0^2 
+ \frac{1}{8}\kappa_0 \kappa_2 + \frac{1}{9} \kappa_1^2 
\) \, dA.
\end{aligned}
\end{equation*}
The residues of $B_M(s)$ for $s=-1$,$-3$,$-5$ were computed by
Brylinski in \cite{bry} (see also \cite{o2018regularized}).
Here, we took the opportunity to correct the coefficient
of $\Res_{s=-5} B_M$ for knots in $\RR^3$.

Therefore
\begin{equation*}
\begin{aligned}
\Res_{s=3} B_1 &= 48 \, \mathrm{Length}(M),\\    
\Res_{s=1} B_1 &= 0, \quad\text{and}\\
\Res_{s=-1} B_1 &= 0.
\end{aligned}
\end{equation*}

Similarly,

\begin{equation*}
\begin{aligned}
\Res_{s=3} B_2 &= 0,\\    
\Res_{s=1} B_2 &= -2\int_M  \kappa_0^2  \, dA, \quad\text{and}\\
\Res_{s=-1} B_2 &= \int_M \(
\frac{3}{4} \kappa_0^4 - \kappa_0^2 \tau_0^2 
+ \kappa_0 \kappa_2 
\) \, dA.
\end{aligned}
\end{equation*}

So the result follows from Equation (\ref{{E:bbcl}}).

\end{proof}

\bibliographystyle{amsplain}
\bibliography{v7-bbcl}

\end{document}